\documentclass{amsart}
\usepackage{graphicx}
\newtheorem{thm}{Theorem}[section]
\newtheorem{cor}[thm]{Corollary}
\newtheorem{lem}[thm]{Lemma}
\newtheorem{prop}[thm]{Proposition}
\theoremstyle{definition}
\newtheorem{defn}[thm]{Definition}
\newtheorem{rem}[thm]{Remark}

\numberwithin{equation}{section}

\begin{document}

\title[]{Local Hardy-Littlewood maximal operator in variable Lebesgue spaces}%
\author{Amiran Gogatishvili, Ana Danelia, Tengiz Kopaliani}%

\address{Institute of Mathematics of the Academy of Sciences of the
Czech Republic, Zitna 25, 11567 Praha 1, Czech Republic}
\email{gogatish@math.cas.cz}

\address{Faculty of Exact and Natural Sciences,
 Tbilisi State University,
 Chavchavadze St.1, Tbilisi 0128 Georgia}%
\email{tengiz.kopaliani@tsu.ge}%

\address{Faculty of Exact and Natural Sciences,
 Tbilisi State University,
 Chavchavadze St.1, Tbilisi 0128 Georgia}%
\email{ana.danelia@tsu.ge}%

\thanks{The research was supported   by Shota Rustaveli National Science Foundation grant
no.13/06 (Geometry of function spaces, interpolation and embedding
theorems). The research  of the first author was  partially supported by the grant
201/08/0383 and 13-14743S of the Grant Agency of the Czech Republic and RVO:
67985840.}%

\begin{abstract}
We investigate the class $\mathcal{B}^{loc}(\mathbb{R}^{n})$ of
exponents $p(\cdot)$ for which the local Hardy-Littlewood maximal
operator is bounded in variable exponent Lebesgue spaces
$L^{p(\cdot)}(\mathbb{R}^{n})$. Littlewood-Paley square function
characterization  of $L^{p(\cdot)}(\mathbb{R}^{n})$ spaces with the
above class of exponent are also obtained.

\end{abstract}
\maketitle
\section{Introduction}

The variable exponent Lebesgue spaces $L^{p(\cdot)}(\mathbb{R}^{n})$   and the corresponding variable exponent Sobolev spaces $W^{k,p(\cdot)}$ are of interest for their applications to the problems in  fluid dynamics \cite{Ru1, Ru2},  partial differential equations with non-standard growth condition and calculus of variations  \cite{AM1, AM2, F, FMP}, image processing \cite{CLR, HHLT, LLP}.

The boundedness of Hardy-Littlewood maximal operator is very important tool to get boundedness of more complicated operators such as singular integral operators, commutators of singular integrals, Riesz potential and many another operators.
Conditions for the boundedness of the Hardy-Littlewood   maximal
operator on  variable exponent Lebesgue spaces
$L^{p(\cdot)}(\mathbb{R}^{n})$ have been studied in \cite{D1, D2,
CFN, N, K2, KK, L}. For an overview we refer to the monograph
\cite{DHHR}.

Let $p:\mathbb{R}^{n}\longrightarrow[1,\infty)$ be a measurable
function. Denote by $L^{p(\cdot)}(\mathbb{R}^{n})$ the space of all
measurable functions $f$ on $\mathbb{R}^{n}$ such that for some
$\lambda>0$
$$\int_{\mathbb{R}^{n}}\left|\frac{f(x)}{\lambda}
\right|^{p(x)}dx<\infty,$$ with the norm
$$\left\|f\right\|_{p(\cdot)}=
\inf\left\{\lambda>0:\int_{\mathbb{R}^{n}}
\left|\frac{f(x)}{\lambda}\right|^{p(x)}dx\leq 1\right\}.$$

Given a locally integrable function $f$ on $\mathbb{R}^{n},$ the
Hardy-Littlewood maximal operator $M$ is defined by the equality
$$ Mf(x)=\sup\frac{1}{|Q|}\int_{Q}|f(y)|dy,$$
where the supremum is taken over all cubes $Q$ containing $x.$
Throughout the paper, all cubes are assumed to have their sides
parallel to the coordinate axes.

Let $f$ be locally integrable function $f$ on $\mathbb{R}^{n}.$  We
consider the local variant of the Hardy-Littlewood maximal operator given
by
$$
M^{loc}f(x)=\sup\limits_{Q\ni
x,|Q|\leq1}\frac{1}{|Q|}\int_{Q}|f(y)|dy.
$$

Denote by $\mathcal{B}(\mathbb{R}^{n})$
($\mathcal{B}^{loc}(\mathbb{R}^{n})$) the class of all measurable
functions 
$p:\mathbb{R}^{n}\longrightarrow[1,\infty)$ for which operator $M$ (operator $M^{loc}$) is
bounded on $L^{p(\cdot)}\left(\mathbb{R}^{n}\right).$ Given any
measurable function $p(\cdot),$ let $p_{-}=\inf_{x\in\mathbb{R}^{n}}
p(x)$ and $p_{+}=\sup_{x\in\mathbb{R}^{n}} p(x).$  Below we assume
that $1<p_{-}\leq p_{+}<\infty.$

 It has been proved by
Diening \cite{D1}
 that if $p(\cdot)$ satisfies the following uniform continuity condition
 \begin{equation}  \label{1.1}
 |p(x)-p(y)|\leq\frac{c}{\log(1/|x-y|)},\,\,\,|x-y|<1/2,
 \end{equation}
 and if $p(\cdot)$ is a constant outside some large ball, then
 $p(\cdot)\in\mathcal{B}(\mathbb{R}^{n}).$ After that the second condition
 on $p(\cdot)$ has been improved independently by Cruz-Uribe, Fiorenza,
 and Neugebauer \cite{CFN} and Nekvinda \cite{N}. It is shown in \cite{CFN} that if $p(\cdot)$
 satisfies \eqref{1.1} and
 \begin{equation} \label{1.2}
 |p(x)-p_{\infty}|\leq\frac{c}{\log(e+|x|)}
 \end{equation}
 for some $p_{\infty}>1,$ then  $p(\cdot)\in\mathcal{B}(\mathbb{R}^{n}).$
 In \cite{N}, the boundedness of $M$ is deduced from  (1.1)  and the integral
 condition more
 general than \eqref{1.2} condition: there exist constants
 $c,p_{\infty},$  such that $0<c<1,p_{\infty}>1,$ and
 $$\int_{\mathbb{R}^{n}}c^{\frac{1}{|p(x)-p_{\infty}|}}dx<\infty.$$

The condition \eqref{1.1} ia named the local log-H\"older continuity
condition and the condition \eqref{1.2} the log-H\"older decay condition (at infinity).
The conditions \eqref{1.1} and \eqref{1.2} together  are named global log-H\"older continuity condition.
this conditions are connected to the geometry of the space $L^{p(\cdot)}(\mathbb{R}^{n})$.

By $\mathcal{X}^{n}$ we denote the set of all open cubes in
$\mathbb{R}^{n}$ and by $\mathcal{Y}^{n}$ ($\mathcal{Y}^{n}_{loc}$)
we denote the set of all families $\mathcal{Q}=\{Q_{i}\}$  of
disjoint, open cubes in $\mathbb{R}^{n}$ (with measure less than
$1$) such that $\bigcup Q_i=\mathbb{R}^{n}$ .

Everywhere below by  $l_{\mathcal{Q}}$ we denote a Banach sequential
space (BSS).   Let $\{e_{Q}\}$ be standard unit vectors in
$l_{\mathcal{Q}}.$

\begin{defn}
Let $l=\{l_{\mathcal{Q}}\}_{\mathcal{Q}\in\mathcal{Y}^{n}}$
($l=\{l_{\mathcal{Q}}\}_{\mathcal{Q}\in\mathcal{Y}^{n}_{loc}}$) be a
family of BSSs. A space  $L^{p(\cdot)}(\mathbb{R}^{n})$  is said to
satisfy a uniformly upper (lower) $l-$estimate ($l_{loc}-$estimate )
if there exists a constant $C>0$ such that for every $f\in
L^{p(\cdot)}(\mathbb{R}^{n})$ and $\mathcal{Q}\in \mathcal{Y}^{n}$
($\mathcal{Q}\in \mathcal{Y}^{n}_{loc}$) we have
$$
\|f\|_{p(\cdot)}\leq
C\|\sum_{Q_{i}\in\mathcal{Q}}\|f\chi_{Q_{i}}\|_{p(\cdot)}\cdot
e_{Q_{i}}\|_{l_{\mathcal{Q}}}\,\,
\left(\|\sum_{Q_{i}\in\mathcal{Q}}\|f\chi_{Q_{i}}\|_{p(\cdot)}\cdot
e_{Q_{i}}\|_{l_{\mathcal{Q}}}\leq C\|f\|_{p(\cdot)}\right).
$$
\end{defn}

Definition 1.1 was introduced by Kopaliani in \cite{K1}. The idea of definition 1.1 is
simply to generalize the following property of the Lebesgue-norm:
$$
\|f\|_{L^{p}}^{p}=\sum_{i}\|f\chi_{\Omega_{i}}\|_{L^{p}}^{p}
$$
for a partition of $\mathbb{R}^{n}$ into measurable sets
$\Omega_{i}.$

 Let $p(\cdot)\in\mathcal{B}(\mathbb{R}^{n}).$ For any $\mathcal{Q}\in \mathcal{Y}^{n}$
 we define the space $l^{\mathcal{Q},p(\cdot)}$ by
$$
l^{\mathcal{Q},p(\cdot)}:=\left\{\overline{t}=\{t_{Q}\}_{Q\in\mathcal{Q}}:\,\,
\sum_{Q\in\mathcal{Q}}|t_{Q}|^{p_{Q}}<\infty\right\},
$$
equipped with the Luxemburg's norm, where the numbers $p_{Q}$ are
defined as  $\frac{1}{p_{Q}}=\frac{1}{|Q|}\int_{Q}\frac{1}{p(x)}dx$.
Analogously we define the space $l^{\mathcal{Q},p'(\cdot)}$ where
$\frac{1}{p(t)}+\frac{1}{p'(t)}=1,\,\,\,t\in\mathbb{R}^{n}.$

Note that if $p(\cdot)\in\mathcal{B}(\mathbb{R}^{n})$ then for
simple functions we have uniformly lower and upper
$l=\{l^{\mathcal{Q},p(\cdot)}\}_{\mathcal{Q}\in\mathcal{Y}^{n}}$
estimates.

\begin{thm}\label{thm1.2}
Let $p(\cdot)\in\mathcal{B}(\mathbb{R}^{n})$ then uniformly
\begin{equation} \label{1.3}
\|\sum_{Q\in\mathcal{Q}}t_{Q}\chi_{Q}\|_{p(\cdot)}\asymp
\|\sum_{Q\in\mathcal{Q}}t_{Q}\|\chi_{Q}\|_{p(\cdot)}e_{Q}\|_{l^{\mathcal{Q},p(\cdot)}}
\end{equation}
and
\begin{equation} \label{1.4}
\|\sum_{Q\in\mathcal{Q}}t_{Q}\chi_{Q}\|_{p'(\cdot)}\asymp
\|\sum_{Q\in\mathcal{Q}}t_{Q}\|\chi_{Q}\|_{p'(\cdot)}e_{Q}\|_{l^{\mathcal{Q},p'(\cdot)}}.
\end{equation}
\end{thm}

 Above  theorem  is another version of necessary part of Diening's Theorem 4.2
 in \cite{D2} (proof may be found in \cite{K3}). Note that
conditions(1.3) and (1.4)  in general do not imply
$p(\cdot)\in\mathcal{B}(\mathbb{R}^{n}).$ The proof (see in  \cite{K5})
relies on the example constructed by Lerner in \cite{L}. We give the proof of this fact also here.

Let $E=\cup_{k\geq1}(e^{k^{3}},e^{k^{3}e^{1/k^2}})$ and
 \begin{equation}
 p_{0}(x)=\int_{|x|}^{\infty}\frac{1}{t\log t}\chi_{E}(t)dt.
\end{equation}
 There exist $\alpha>1$ and $\beta_{0}\,(1/\alpha<\beta_{0}<1)$ such that $p_{0}(\cdot)+\alpha\in\mathcal{B}(\mathbb{R})$ and
$\beta_{0}(p_{0}(\cdot)+\alpha)\notin \mathcal{B}(\mathbb{R})$ (see
\cite[Theorem 1.7]{L}). Note that  for a space
$L^{p(\cdot)}(\mathbb{R})$, with $p(\cdot)=p_{0}(\cdot)+\alpha$
there exists a family
$l=\{l_{\mathcal{Q}}\}_{\mathcal{Q}\in\mathcal{Y}^{n}}$
 of BSSs for which $L^{p(\cdot)}(\mathbb{R})$ satisfies uniformly lower and upper
 $l-$estimate (see \cite[Proposition~3.2]{K5}). From \eqref{1.3} we have $l_{\mathcal{Q}}\cong l^{\mathcal{Q},p(\cdot)}$ and consequently
 we have
\begin{equation}\label{1.6}
\|f\|_{p(\cdot)}\asymp
\|\sum_{Q\in\mathcal{Q}}\|f\chi_{Q}\|_{p(\cdot)}e_Q\|_{l^{\mathcal{Q},p(\cdot)}}.
\end{equation}
Note that for all $1>\beta>\frac1{p_{-}}$
\begin{equation}\label{1.7}
\|f^{\frac1{\beta}}\|_{\beta p(\cdot)}^{\beta}=\|f\|_{p(\cdot)}\\
\end{equation}
and
\begin{equation} \label{1.8}
\|\{t_Q\}\|_{l^{\mathcal{Q},p(\cdot)}}=\left\|\left\||t_Q|^{\frac{1}{\beta}}\right\}\right\|_{l^{\mathcal{Q},p(\cdot)}}^{\beta}.
\end{equation}
From \eqref{1.6}, \eqref{1.7} and \eqref{1.8} we have
$$ \|g\|_{\beta p(\cdot)}\asymp
\|\sum_{Q\in\mathcal{Q}}\|g\chi_{Q}\|_{p(\cdot)}e_Q\|_{l^{\mathcal{Q},\beta p(\cdot)}}.
$$
 for $g\in L^{\beta p(\cdot)}(\mathbb{R})$ and the space $ L^{\beta p(\cdot)}(\mathbb{R})$
 satisfies uniformly lower and upper  $l^\beta$-estimates, where $l_Q^\beta=l^{Q,\beta p(\cdot)}$.

 Note that $\frac1{(\beta p(\cdot))_Q}+\frac1{((\beta p(\cdot))')_Q}=1$ and
 $\left(l^{\mathcal{Q},\beta p(\cdot)}\right)'=l^{\mathcal{Q},(\beta p(\cdot))'}$.
 Thus the space $ \left(L^{\beta p(\cdot)}(\mathbb{R})\right)'$  satisfies uniformly lower and
 upper $(l^\beta)'$-estimates, where $(l^\beta)_Q'=l^{Q,(\beta p(\cdot))'}$ and \eqref{1.3} and \eqref{1.4}
 are valid for any $\beta p(\cdot)$, $(\beta p(\cdot))'$, where $1>\beta>\frac1{p_{-}}$. Consequently for
  exponent $\beta_0p(\cdot)$ \eqref{1.3} and \eqref{1.4} are valid but, $\beta_0p(\cdot)\notin \mathcal{B}(\mathbb{R})$.

\begin{rem} \label{rem1.1} Let $p(\cdot)$ be  global  log-H\"{o}lder continuous
function. 
Then there exists family
$l=\{l_{\mathcal{Q}}\}_{\mathcal{Q}\in\mathcal{Y}^{n}}$
 of BSSs for which $L^{p(\cdot)}(\mathbb{R}^{n})$ satisfies uniformly lower and upper
 $l-$estimates (see \cite[Proposition~3.4]{K5}). As we already mentioned it was show
 in \cite{CFN} that $p(\cdot)\in\mathcal{B}(\mathbb{R}^{n})$  and by Theorem~\ref{thm1.2} \eqref{1.3}
 holds and therefore we have  $l_{\mathcal{Q}}\cong l^{\mathcal{Q},p(\cdot)}$  and consequently
  \begin{equation}\label{1.9}
\|f\|_{p(\cdot)}\asymp
\|\sum_{Q\in\mathcal{Q}}\|f\chi_{Q}\|_{p(\cdot)}e_Q\|_{l^{\mathcal{Q},p(\cdot)}}.
\end{equation}
\end{rem}

\begin{rem} \label{rem1.2} Let $\mathcal{Q}=\{Q_{i}\}$ be a partition of
$\mathbb{R}^{n}$ into equal sizes cubes, ordered so that $i>j$ if
$\mbox{dist}(0,Q_{i})>\mbox{dist}(0,Q_{j}).$  Let $p(\cdot)$ be
global log-H\"{o}lder continuous. Then
\begin{equation}\label{1.10}
\|f\|_{p(\cdot)}\approx\left(\sum_{i}\|f\chi_{Q_{i}}\|_{p(\cdot)}^{p_{\infty}}\right)^{1/p_{\infty}}.
\end{equation}
 This was shown in \cite[Theorem~2.4]{H}. This statement also follows from Remark~\ref{rem1.1}.
  Indeed, if we have a partition   $\mathcal{Q}=\{Q_{i}\}$ with equal sizes cubes and it is ordered as above
  by using \cite[Theorem~4.3]{N1} we can show that $l^{p_\infty}\cong l^{\mathcal{Q},p(\cdot)}$
  and consequently from \eqref{1.9} we get\eqref{1.10}.
\end{rem}

 By $\mathcal{AC}$  we denote the set of  exponents
$p:\mathbb{R}\rightarrow[1,+\infty)$  of the form
$p(x)=p+\int_{-\infty}^{x}l(u)du,$ where
$\int_{-\infty}^{+\infty}|l(u)|du<+\infty.$

Note that example of exponent constructed by Lerner and mentioned
above belongs to class $\mathcal{AC}.$ In general we have the
following

\begin{prop}\cite[Proposition~3.2]{K5} Let  $p(\cdot)\in\mathcal{AC}.$ Then exists
family $l=\{l_{\mathcal{Q}}\}_{\mathcal{Q}\in\mathcal{Y}^{n}}$
 of BSSs for which $L^{p(\cdot)}(\mathbb{R}^{n})$ satisfies uniformly lower and upper
 $l-$estimate.
\end{prop}

In many applications it is enough to study only boundadnes of local Hardy-Littlewood maximal
operator rather the  Hardy-Littlewood maximal
operator. For example in the Littlewood-Paley theory we need local  Hardy-Littlewood maximal
operator. In the weighted Lebesgue spaces behavior of local  Hardy-Littlewood maximal
operator was studied by Rychkov in \cite{Ry}.

In this paper we investigate the class $\mathcal{B}^{loc}(\mathbb{R}^{n})$ of
exponents $p(\cdot)$ for which the local Hardy-Littlewood maximal
operator is bounded in variable exponent Lebesgue space $L^{p(\cdot)}(\mathbb{R}^{n})$.
 Using the obtained results we give Littlewood-Paley square-function characterization  of
 the variable exponent Lebesgue spaces $L^{p(\cdot)}(\mathbb{R}^{n})$ with the above class of
exponent.

The paper is organized as follows. In Section 2 we give main results. In section 3  we give application in
the  Littlewood-Paley theory and in last section we give outlines of the proof  of the Theorem~2.2
which is local version of the Dieninges theorem from \cite{D2}.

\section{Main results}

For any family of pairwise disjoint cubes $\mathcal{Q}$ and $f\in
L^{1}_{loc}$ we define the averaging operator
$$T_{\mathcal{Q}}f=\sum_{Q\in\mathcal{Q}}\chi_{Q}M_{Q}f$$
where $M_{Q}f=|Q|^{-1}\int_{Q}f(x)dx.$

 We say that exponent $p(\cdot)$ is of class $\mathcal{A}$ ( class $\mathcal{A}^{loc}$) if and
 only if there exists $C>0$ such that for all $\mathcal{Q}\in
 \mathcal{Y}^{n}$  ($\mathcal{Q}\in
 \mathcal{Y}^{n}_{loc}$) and all $f\in L^{p(\cdot)}(\mathbb{R}^{n})$
 $$
 \|T_{\mathcal{Q}}\|_{p(\cdot)}\leq
 C\|f\|_{p(\cdot)},
 $$
 i.e. the averaging operators  $T_{\mathcal{Q}}$ are uniformly
 continuous on $ L^{p(\cdot)}(\mathbb{R}^{n}).$

A necessary and sufficient condition on  $p(\cdot)$ for which
operator $M$ is bounded in $L^{p(\cdot)}(\mathbb{R}^{n})$ is given
by Diening in \cite{D2}. It states that
$p(\cdot)\in\mathcal{B}(\mathbb{R}^{n})$ if the averaging operators
$T_{\mathcal{Q}}$ are uniformly
 continuous on $ L^{p(\cdot)}(\mathbb{R}^{n})$
 with respect to all families $\mathcal{Q}$ of disjoint cubes.
 This concept provides the following characterization of when the
 maximal operator is bounded.
\begin{thm} ( \cite[Theorem 8.1]{D2}). Let $1<p_{-}\leq
 p_{+}<\infty.$ The following are equivalent:

 1) $p(\cdot)$ is of class $\mathcal{A};$

 2) $M$ is bounded on  $ L^{p(\cdot)}(\mathbb{R}^{n});$

 3) $(M(|f|^{q}))^{1/q}$ is bounded on $
 L^{p(\cdot)}(\mathbb{R}^{n})$ for some $q>1,$ ("left-openness");

 4)$M$ is bounded on $
 L^{p(\cdot)/q}(\mathbb{R}^{n})$ for some $q>1,$ ("left-openness");

 5) $M$ is bounded on $L^{p'(\cdot)}(\mathbb{R}^{n}).$

\end{thm}

The following theorem is analog of Dienings theorem for local
maximal function $M^{loc}.$  The proof is principally the same as of
Theorem 2.1 with minor modification  of  some facts. The convenience
for  reader in appendix we sketch of proof.

\begin{thm}   Let $1<p_{-}\leq
 p_{+}<\infty.$ The following are equivalent:

 1) $p(\cdot)$ is of class $\mathcal{A}^{loc};$

 2) $M^{loc}$ is bounded on  $ L^{p(\cdot)}(\mathbb{R}^{n});$

 3) $(M^{loc}(|f|^{q}))^{1/q}$ is bounded on $
 L^{p(\cdot)}(\mathbb{R}^{n})$ for some $q>1,$ ("left-openness");

 4)$M^{loc}$ is bounded on $
 L^{p(\cdot)/q}(\mathbb{R}^{n})$ for some $q>1,$ ("left-openness");

 5) $M^{loc}$ is bounded on $L^{p'(\cdot)}(\mathbb{R}^{n}).$
\end{thm}

We say that $dx$ satisfies the condition $A_{p(\cdot)}$  (condition
$A_{p(\cdot)}^{loc}$ ) if there exists $C>0$ such that for any cube
$Q$ (for any cube $Q$ with $|Q|\leq1$)
$$
\frac{1}{|Q|}\|\chi_{Q}\|_{p(\cdot)}\|\chi_{Q}\|_{p'(\cdot)}\leq C.
$$

Using Theorem 2.1-2.2  we obtain some subclass of
$\mathcal{B}(\mathbb{R}^{n})$ and
$\mathcal{B}^{loc}(\mathbb{R}^{n}).$

\begin{thm}
Let $1<p_{-}\leq p_{+}<\infty$ and there exists family
$l=\{l_{\mathcal{Q}}\}_{\mathcal{Q}\in\mathcal{Y}^{n}}$ (family
$l=\{l_{\mathcal{Q}}\}_{\mathcal{Q}\in\mathcal{Y}^{n}_{loc}}$)
 of BSSs for which $L^{p(\cdot)}(\mathbb{R}^{n})$ satisfies uniformly lower and upper
 $l-$estimate ($l_{loc}-$ estimate). Then operator $M$ (operator $M^{loc}$)
 is bounded in  $L^{p(\cdot)}(\mathbb{R}^{n})$ if and only if $dx\in
 A_{p(\cdot)}$ ($dx\in
 A_{p(\cdot)}^{loc}$).
\end{thm}

Proof. The proof for  operator $M^{loc}$ is  the same as for
operator $M.$ Let $dx\in
 A_{p(\cdot)}^{loc}$. Using H\"{o}lders inequality we get
 $$
\frac{1}{|Q|}\int_{Q}|f(x)|dx\leq
C\frac{\|f\chi_{Q}\|_{p(\cdot)}}{\|\chi_{Q}\|_{p(\cdot)}}.
$$

For $\mathcal{Q}\in\mathcal{Y}^{n}_{loc}$ and $f\in
L^{p(\cdot)}(\mathbb{R}^{n})$ we have
$$
\left\|\sum_{Q\in\mathcal{Q}}\chi_{Q}\frac{1}{Q}\int_{Q}f(x)dx\right\|_{p(\cdot)}\leq
\|\sum_{Q\in\mathcal{Q}}\|f\chi_{Q}\|e_{Q}\|_{l_{\mathcal{Q}}}\leq
 C\|f\|_{p(\cdot)}.
 $$

 The necessary part of theorem is obvious.

\begin{thm}
$\mathcal{B}(\mathbb{R}^{n})\neq \mathcal{B}^{loc}(\mathbb{R}^{n})$
\end{thm}

Proof. Let us consider the exponent
$p(\cdot)=\beta_{0}(p_{0}(\cdot)+\alpha)$ where $p_{0}(\cdot)$ is
defined by (1.5). Let fix $\alpha>1$ and
$\beta_{0}\,(1/\alpha<\beta_{0}<1)$ such that the exponent
$p(\cdot)$ does not belong to the class $\mathcal{B}(\mathbb{R}).$
Since $p(\cdot)\in\mathcal{AC}$ we can conclude that for
$L^{p(\cdot)}(\mathbb{R})$ there exists family
$l=\{l_{\mathcal{Q}}\}_{\mathcal{Q}\in\mathcal{Y}^{n}}$
 of BSSs for which $L^{p(\cdot)}(\mathbb{R})$ satisfies uniformly lower and upper
 $l-$estimate. For $p(\cdot)$ the condition (1.1) is fulfilled, so
 it is easy to show that for this exponent $p(\cdot)$ condition $A_{p(\cdot)}^{loc}$
  is satisfied.Therefore by Theorem 2.2  $p(\cdot)\in\mathcal{B}^{loc}(\mathbb{R}).$

  In the class $\mathcal{B}^{loc}(\mathbb{R})$ there
 exist exponents that have arbitrary slow decreasing order in
 infinity.  To show this fact we rely on the simple observation.
 Indeed, let for $L^{p(\cdot)}(\mathbb{R})$ there exists family
 $l=\{l_{\mathcal{Q}}\}_{\mathcal{Q}\in\mathcal{Y}^{n}}$
 of BSSs for which $L^{p(\cdot)}(\mathbb{R}^{n})$ satisfies uniformly lower and upper
 $l-$ estimate and $\omega:\mathbb{R}\rightarrow\mathbb{R};\,\omega(-\infty)=-\infty,\,\omega(+\infty)=+\infty$ is
 strictly increasing absolutely continuous mapping. Then there exists family
 $l_{\omega}$
 of BSSs for which $L^{p(\omega(\cdot))}(\mathbb{R})$  satisfies uniformly lower and upper
 $l_{\omega}-$ estimate (see \cite{K5}).

  Consider the exponent from Theorem 2.4. Let $t_{k}=e^{k^{3}},\,\,m_{k}=e^{k^{3}}e^{1/k^{2}},\,k\geq1.$
   Let us construct new points $t'_{k},\,m'_{k},\,k\geq1$ so that
  $m'_{k}-t'_{k}=m_{k}-t_{k}$  and $t'_{k+1}>m'_{k}.$
  Let us now construct the pairwise linear continuous function
  $\omega$ in the following way: $\omega(x)=x$ if $x\leq0,$ $\omega(t_{k})=t'_{k},\,\,\omega(m_{k})=m'_{k};\,k\geq1.$
   We can choose the points $t'_{k},\,m'_{k}$  so that
   $(m'_{k+1}-t'_{k})/(m_{k+1}-t_{k})$ was arbitrary large. Note
  that exponents  $p(w^{-1}(\cdot))$ and $p(\cdot)$ has the same
  local behavior   but the decreasing order in infinity of
   $p(w^{-1}(\cdot))$ is very slow.

   Let now consider the case $n\geq2.$ Let
   $D=\cup_{k=1}^{\infty}[2k-1,2k]\times[0,1]^{n-1}.$ Consider
   non-trivial exponent $p(\cdot)$ that satisfies global log-H\"{o}lder
   condition and is constant on the set $\mathbb{R}^{n}\backslash D.$

   Let $\{m_{k}\}$ be the strictly increasing sequence of integers.
   Consider the bijection
   $\omega:\mathbb{R}^{n}\rightarrow\mathbb{R}^{n}$ that for each
   $k\in\mathbb{N}$   has the
   form $\omega(x)=x-(m_{k},0,...,0)$  on the set $[2k-1,2k]\times[0,1]^{n-1}.$ We can choose the
   sequence $\{m_{k}\}$ so
   that $p(\omega(\cdot))\overline{\in}\mathcal{B}(\mathbb{R}^{n})$  but
   $p(\omega(\cdot))\in\mathcal{B}^{loc}(\mathbb{R}^{n}).$
$\Box$

 Note that only the condition $dx\in A_{p(\cdot)}^{loc}$ (even $dx\in A_{p(\cdot)}$ ) does not
 guarantee in general
 $p(\cdot)\in\mathcal{B}^{loc}(\mathbb{R}^{n}).$
 The corresponding example see in \cite{K4}.

\section{Some applications}

In this section, we give Littlewood-Paley square-function
characterization of $L^{p(\cdot)}(\mathbb{R}^{n})$ when
$p(\cdot)\in\mathcal{B}^{loc}(\mathbb{R}^{n}).$
 Let us recall the definition of local  Muckenhoupt weights.
The weight class $A_{p}^{loc}\,\,(1<p<\infty)$ to consists of all
nonnegative locally integrable functions $w$ on $\mathbb{R}^{n}$ for
which
$$
A_{p}^{loc}(w)=\sup\limits_{|Q|\leq1}\frac{1}{|Q|^{p}}\int_{Q}w(x)dx\left(w(x)^{-p'/p}dx\right)^{p/p'}<\infty.
$$

Extending the suprema from $|Q|\leq1$ to all $Q$ gives the
definition  of the usual  classes  $A_{p}$. It follows directly from
definition that $A_{p}\subset A_{p}^{loc}.$  The littlwood-Paley
theory for weight lebesgue space $L^{p}_{w}$ with local  Muckenhoupt
weights was investigate by Rychkov in \cite{Ry}.  For more details
for $A_{p}^{loc}$ weights we refer paper \cite{Ry}.

 Below we formulate analog of Rubio de
Francia theorem for variable exponent case. Hereafter, $\mathcal{F}$
will denote a family of ordered pairs of non-negative, measurable
functions $(f,g).$ If we say that for some $p$, $1<p<\infty,$ and
$w\in A^{loc}_{p}$

\begin{equation}
\int_{\mathbb{R}^{n}}f(x)^{p}w(x)dx\leq
C\int_{\mathbb{R}^{n}}g(x)^{p}w(x)dx,\quad (f,g)\in\mathcal{F},
\end{equation}
we mean that this inequality holds for any $(f,g)\in\mathcal{F}$
such that the left-hand side is finite, and that the constant $C$
depends only on $p$ and the constant $A^{loc}_{p}(w).$

\begin{thm}.
Given a family $\mathcal{F},$ assume that (3.1) holds for some
$1<p_{0}<\infty,$ for every weight $\omega\in A_{p_{0}}^{loc}$ and
 for all $(f,g)\in \mathcal{F}$. Let
$p(\cdot)$ be such that there
 exists $1<p_{1}<p_{-},$ with $(p(\cdot)/p_{1})'\in \mathcal{B}^{loc}(\mathbb{R}^{n}).$
Then
$$
\|f\|_{p(t)}\leq C\|g\|_{p(t)}
$$
for all $(f,g)\in\mathcal{F}$ such that $f\in
L^{p(t)}(\mathbb{R}^{n}).$ Furthermore, for every $0<q<\infty$ and
sequence $\{(f_{j},g_{j})\}_{j}\subset\mathcal{F},$
$$
\left\|\left(\sum_{j}(f_{j})^{q}\right)^{1/q}\right\|_{p(t)}\leq
C\left\|\left(\sum_{j}(g_{j})^{q}\right)^{1/q}\right\|_{p(t)}.
$$
\end{thm}

 In case when $w\in A_{p_{0}}$ and
$(p(\cdot)/p_{1})'\in\mathcal{B}(\mathbb{R}^{n})$ Theorem 3.1 proved
in \cite{CFMP} (Theorem 1.3) (see also proof Theorem 3.25 in
\cite{CMP}).  Note that the collection of all cubes $Q$ with
$|Q|\leq1$ form the Muckenhoupt basis , that is for each
$p,\,1<p<\infty,$ and for every $w\in A_{p}^{loc},$ the maximal
operator $M_{loc}$ is bounded on $L^{p}_{w}(\mathbb{R}^{n})$
(\cite{Ry}, lemma 2.11). The theorem 3.1 follows from theorem 2.2
and   extrapolation theorem for general Banach function spaces
(\cite{CMP}, Theorem 3.5).

We give a number of applications of Theorem 3.1. It is well known
(see \cite{Ry})  that for $1<p<\infty$ and for $w\in A_{p}^{loc},$
$$
\int_{\mathbb{R}^{n}}M^{loc}f(x)^{p}w(x)dx\leq C
\int_{\mathbb{R}^{n}}f(x)^{p}w(x)dx.
$$

From Theorem 3.1 with the pairs $(M^{loc}f,|f|),$ we get
vector-valued inequalities for $M^{loc}$  on
$L^{p(\cdot)}(\mathbb{R}^{n}),$ provided there exists
$1<p_{1}<p_{-}$ with $(p(\cdot)/p_{1})'\in
\mathcal{B}^{loc}(\mathbb{R}^{n});$ by Theorem 2.2, this is
equivalent to $p(\cdot)\in\mathcal{B}^{loc}(\mathbb{R}^{n}).$ We
obtain following local version of the Fefferman-Stein vector-valued
maximal theorem:
\begin{cor} Let $p(\cdot)\in \mathcal{B}^{loc}(\mathbb{R}^{n}).$
Then for all $1<q<\infty,$
$$
\left\|\left(\sum_{j}(M^{loc}f_{j})^{q}\right)^{1/q}\right\|_{p(t)}\leq
C\left\|\left(\sum_{j}(g_{j})^{q}\right)^{1/q}\right\|_{p(t)}.
$$
\end{cor}

Let $1<p<\infty$ and $w\in A_{p}^{loc}.$  Let $\varphi_{0}\in
C_{0}^{\infty}$ have nonzero integral, and
$\varphi(x)=\varphi_{0}(x)-2^{-n}\varphi_{0}(\frac{x}{2}),\,\,x\in\mathbb{R}^{n}.$
Consider the square operator $S=S_{\varphi_{0},\varphi}$ given by
\begin{equation}
S(f)=\left( \sum_{j=0}^{+\infty}|\varphi_{j}\ast
f|^{2}\right)^{1/2}\,\,\,\,(f\in L^{p}_{w}(\mathbb{R}^{n})),
\end{equation}
where $\varphi_{j}(x)=2^{jn}\varphi(2^{j}x),\,\,j\in\mathbb{N}.$
 Then
$$
\|S(f)\|_{L^{p}_{w}}\approx\|f\|_{L^{p}_{w}},\,\,\,\mbox{all}\,\,\,f\in
L^{p}_{w}(\mathbb{R}^{n}).
$$
 (For details, see \cite{Ry}). Therefore
by theorem 3.1 we have following Littlewood-Paley square-function
characterization of $L^{p(\cdot)}(\mathbb{R}^{n}).$

\begin{cor}  Let $p(\cdot)\in \mathcal{B}^{loc}(\mathbb{R}^{n}).$
 Let
$\varphi_{0}\in C_{0}^{\infty}$ have nonzero integral, and
$\varphi(x)=\varphi_{0}(x)-2^{-n}\varphi_{0}(\frac{x}{2}).$ Consider
the square operator $S=S_{\varphi_{0},\varphi}$ given by equation
(3.2). Then
$$
\|S(f)\|_{p(\cdot)}\approx\|f\|_{p(\cdot)},\,\,\,\mbox{all}\,\,\,f\in
L^{p(\cdot)}(\mathbb{R}^{n}).
$$
\end{cor}

\section{Appendix}

 Let $\varphi(x,t)=t^{p(x)}$
$t\geq0,\,x\in\mathbb{R}^{n},\,1<p_{-}\leq p_{+}<\infty.$    We need
some notations. For $t\geq0,\,s\geq1,$ we define

$$
\varphi(f)(x):\,\mathbb{R}^{n}\rightarrow[0,+\infty)=\mathbb{R}^{\geq0},\,\,\,(\varphi(f))(x)=\varphi(x,\,|f(x)|),
$$
$$
$$
$$
M_{s,Q\varphi}:\,\mathbb{R}^{n}\rightarrow\mathbb{R}^{\geq0},\,\,\,M_{s,Q\varphi}(t)=\left(\frac{1}{|Q|}\int_{Q}(\varphi(x,t))^{s}dx\right)^{1/s}
$$

$$
M_{Q\varphi}:\,\mathbb{R}\rightarrow\mathbb{R}^{\geq0},\,\,\,M_{Q\varphi}(t)=\left(M_{1,Q\varphi}\right)(t).
$$
Analogously we will  use notation for  the complementary function of
$\varphi$   given by
$\varphi^{\ast}(x,t)=(p(x)-1)p(x)^{-p'(x)}t^{p'(x)}.$

Note that for all cube $Q$ functions
$(M_{s,Q\varphi})(t),\,\,(M_{s,Q\varphi^{\ast}})(t)$ are
$N$-functions  and   satisfy uniformly
$\bigtriangleup_{2}$-condition with respect to $Q$ (see \cite{D2},
Lemma 3.4). In addition we mention following properties of functions
defined above (\cite{D2}, Lemma 3.7):
 let
$s\geq1$ and $Q\in \mathcal{X}^{n},$ then for all $f\in
L^{p(\cdot)}(\mathbb{R}^{n})$ there holds
\begin{equation}
(M_{s,Q\varphi^{\ast}})^{\ast}\left(\frac{1}{2}M_{s,Q}f\right)\leq
M_{s,Q}(\varphi(f)).
\end{equation}
Especially, for all $u>0$
\begin{equation}
(M_{s,Q\varphi^{\ast}})^{\ast}\left(\frac{1}{2}u\right)\leq
M_{s,Q\varphi}(u).
\end{equation}
On the other hand for all $t>0$ the function
$f_{t}=\chi_{Q}\varphi^{\ast}(t)/t$ satisfies
\begin{equation}
(M_{s,Q\varphi^{\ast}})^{\ast}(2M_{s,Q}f_{t})\geq
M_{s,Q}(\varphi(f_{t})).
\end{equation}

For $\mathcal{Q}\in\mathcal{Y}^{n}$ we define the space
$l^{|Q|M_{Q\varphi}}(\mathcal{Q})$

$$
l^{|Q|M_{Q\varphi}}(\mathcal{Q})=\left\{\overline{t}=\{t_{Q}\}_{Q\in\mathcal{Q}}:\,\,\sum_{Q\in\mathcal{Q}}|Q|(M_{Q\varphi})(t_{Q})<\infty\right\},
$$
equipped with the norm
$$
\left\|\overline{t}\right\|_{l^{|Q|M_{Q\varphi}}(\mathcal{Q})}=\inf\left\{\lambda>0:\,\sum_{Q\in\mathcal{Q}}|Q|(M_{Q\varphi})(t_{Q}/\lambda)<1\right\}.
$$

Analogously we define the spaces
$l^{|Q|M_{Q\varphi^{\ast}}}(\mathcal{Q}),$
$l^{|Q|M_{s,Q\varphi}}(\mathcal{Q}),$
$l^{|Q|M_{s,Q\varphi^{\ast}}}(\mathcal{Q}).$

\begin{defn} Let
$$
l^{|Q|(M_{Q\varphi^{\ast}})^{\ast}}(\mathcal{Q})\hookrightarrow
l^{|Q|M_{Q\varphi}}(\mathcal{Q})$$
 are uniformly continuous with respect to
 $\mathcal{Q}\in\mathcal{Y}^{n}_{loc}$
 ($\mathcal{Q}\in\mathcal{Y}^{n}$) i.e.
  for all $A_{1}>0$  there
exists $A_{2}>0$  such that for all
$\mathcal{Q}\in\mathcal{Y}^{n}_{loc}$  (all
$\mathcal{Q}\in\mathcal{Y}^{n}$)  and all sequences
$\{t_{Q}\}_{Q\in\mathcal{Q}}$  there holds
$$
\sum_{Q\in\mathcal{Q}}|Q|(M_{Q\varphi^{\ast}})^{\ast}(t_{Q})\leq
A_{1}\,\,\,\Rightarrow\,\,\,\sum_{Q\in\mathcal{Q}}|Q|(M_{Q\varphi})(t_{Q})\leq
A_{2}.$$
 Then  we say that $M_{Q\varphi}$ is locally  dominated
(dominated) by $(M_{Q\varphi^{\ast}})^{\ast}$  and write
$M_{Q\varphi}\preceq (M_{Q\varphi^{\ast}})^{\ast}(loc)$
($M_{Q\varphi}\preceq (M_{Q\varphi^{\ast}})^{\ast}$).
\end{defn}

Analogously we may define uniformly continuous embedding discrete
function spaces defined above with respect to
$\mathcal{Q}\in\mathcal{Y}^{n}_{loc}$
($\mathcal{Q}\in\mathcal{Y}^{n}$). The basic property of domination
($\preceq$) in a "pointwise" sense is described in original paper
\cite{D2}.
 analogous properties of local domination is essentially based on
 the
following general lemma (note that if  $X=\mathcal{X}^{n}(loc)$ and
$Y=\mathcal{Y}^{n}(loc),$ then $X,\,Y$ are admissible for  Lemma 4.2
)

\begin{lem} (\cite{D2}, Lemma 7.1)
Let $X$ be an arbitrary set. Let $Y$ be a subset of the power set of
$X$ such that $M_{1}\subset M_{2}\in Y$ implies $M_{1}\in Y.$ Let
$\psi_{1},\,\psi_{2}\,:X\rightarrow \mathbb{R}^{\geq}.$  If there
exists $A_{1}>0$  and $A_{2},\,A_{3}\geq0$  such that for all $M\in
Y$
$$
\omega\psi_{1}(\omega)\leq
A_{1}\,\,\,\Rightarrow\,\,\,\,\sum_{\omega\in M}\psi_{2}(\omega)\leq
A_{2}\sum_{\omega\in M}\psi_{1}(\omega)+A_{3}
$$
then there exists $b:\,X\rightarrow\mathbb{R}^{\geq}$  such that for
all $\omega\in X$  holds
\begin{equation}
\psi_{1}(\omega)\leq\frac{A_{1}}{4}\,\,\,\,\,\Rightarrow\,\,\,\,\psi_{2}\leq\max\left\{\frac{4A_{3}}{A_{1}},\,2A_{2}\right\}\psi_{1}+b(\omega)
\end{equation}
and
\begin{equation}
\sup_{M\in Y}\sum_{\omega\in M}b(\omega)\leq A_{3}.
\end{equation}
If on the other hand there exist
$b:\,X\rightarrow\mathbb{R}^{\geq},\,\,A_{1}>0,$ and
$A_{2},\,A_{3}\geq0$  such that (4.4) and (4.5) hold, then for all
$M\in Y$
$$
\psi_{1}(\omega)\leq\frac{A_{1}}{4}\,\,\,\,\,\Rightarrow\,\,\,\,\sum_{\omega\in
M}\psi_{2}(\omega)\leq\max\left\{\frac{4A_{3}}{A_{1}},\,2A_{2}\right\}\omega\psi_{1}(\omega)+A_{3}.
$$
\end{lem}

We can now state characterization of classes  $\mathcal{A}^{loc}$
and $\mathcal{A}.$

\begin{thm}
Exponent $p(\cdot)$ is of class $\mathcal{A}^{loc}$ (of class
$\mathcal{A}$) if and only if
$M_{Q\varphi}\preceq(M_{Q\varphi^{\ast}})^{\ast}(loc)$
($M_{Q\varphi}\preceq(M_{Q\varphi^{\ast}})$)
\end{thm}

The proof of above theorem in case $p(\cdot)$ is of class
$\mathcal{A}$ is  based on  properties (4.1)-(4.3) of $M_{Q\varphi}$
and $(M_{Q\varphi^{\ast}})$  and may use analogously arguments in
local variant.

Inspired by the classical Muckenhoupt class $A_{\infty}$ in
\cite{D2} was defined condition $\mathcal{A}_{\infty}.$  The
importance of our considerations is analogous of definition in local
case.

\begin{defn}
We say that exponent  $p(\cdot)$ is of class
$\mathcal{A}_{\infty}^{loc}$ (class $\mathcal{A}_{\infty}$) if for
any $\varepsilon>0$  there exists $\delta>0$  such that the
following holds: if $N\subset\mathbb{R}^{n}$ is measurable and
$\mathcal{Q}\in \mathcal{Y}^{n}_{loc}$ ($\mathcal{Q}\in
\mathcal{Y}^{n}$) such that
$$
|Q\cap N|\geq\varepsilon|Q|\,\,\,\,\,\mbox{for
all}\,\,\,Q\in\mathcal{Q},
$$
then for any sequence $\{t_{Q}\}_{Q\in\mathcal{Q}}$
$$
\delta\left\|\sum_{Q\in\mathcal{Q}}t_{Q}\chi_{Q}\right\|_{p(\cdot)}\leq\left\|\sum_{Q\in\mathcal{Q}}t_{Q}\chi_{Q\cap
N}\right\|_{p(\cdot)}.
$$
\end{defn}

It is not hard to proof that  if exponent $p(\cdot)$ is of class
$\mathcal{A}^{loc}$ then exponent $p(\cdot)$ is in
$\mathcal{A}^{loc}_{\infty}.$

The important property of exponents from class
$\mathcal{A}_{\infty}$ is that  $\mathcal{A}_{\infty}$ implies
$M_{s,Q\varphi}\preceq M_{Q\varphi}$ for some $s>1.$  The proof of
this result is based on the following lemma.
\begin{lem} (\cite{D2} Lemma 5.5)
Let exponent $p(\cdot)$ is of class $A_{\infty}.$ Then there exists
$\delta>0$ and $A\geq1$ such that for all $Q\in \mathcal{Y}^{n},$
all $\{t_{Q}\}_{Q\in\mathcal{Q}}, t_{Q}\geq 0,$ and all $f\in
L_{loc}^{1}$ with $M_{Q}f\neq0,\,\,Q\in\mathcal{Q},$  holds
$$
\left\|\sum_{Q\in\mathcal{Q}}t_{Q}\left|\frac{f}{M_{Q}f}\right|^{\delta}\chi_{Q}\right\|_{p(\cdot)}\leq
A \left\|\sum_{Q\in\mathcal{Q}}t_{Q}\chi_{Q}\right\|_{p(\cdot)}.
$$
\end{lem}

Note that a very similar argument can be used to obtain  local
version of Lemma 4.5. In original proof of Lemma 4.5 is used
$Q$-dyadic ($Q\in\mathcal{X}^{n}$) maximal function
$M^{\bigtriangleup,Q}$ (\cite{D2}, Definition 5.4). Note that in
fact in proof of Lemma 4.5 it is used local $Q$-dyadic maximal
function, where the supremum  is taken over all $Q$-dyadic cube $Q'$
containing $x$ and $|Q'|\leq|Q|.$ As a consequence of local variant
of Lemma 4.5 we obtain a kind reverse H\"{o}lder estimate for
exponents from class $\mathcal{A}^{loc}.$
\begin{thm}
Let $p(\cdot)\in \mathcal{A}^{loc}.$ Then there exists $s>1,$ such
that $M_{s,Q\varphi}\preceq M_{Q\varphi}(loc).$
\end{thm}

From Theorem 4.3 and Theorem 4.6 for local variant we obtain
\begin{thm}
The following conditions are equivalent

$(a)$  $p(\cdot)$ is of class $\mathcal{A}^{loc}$.

$(b)$  $M_{Q\varphi}\preceq (M_{Q\varphi^{\ast}})^{\ast}(loc)$

$(c)$  There exists $s>1,$ such that $M_{s,Q\varphi}\preceq
M_{Q\varphi}\preceq
(M_{Q\varphi^{\ast}})^{\ast}\preceq(M_{s,Q\varphi^{\ast}})^{\ast}(loc).$
\end{thm}

The key lemma from which was derived original Theorem 3.1 is Lemma
8.7 from \cite{D2}. We formulate analogous statement for local
variant.

\begin{lem}
Let $p(\cdot)\in\mathcal{A}^{loc}.$ Then there exists $s>1$ such
that for all $A_{1}>0$ there exist $A_{2}>0$ such that the following
holds:

For all families
$\mathcal{Q}_{\lambda}\in\mathcal{Y}^{n}_{loc},\,\,\lambda>0,$  with
$$
\sum_{Q\in\mathcal{Q}_{\lambda}}|Q|(M_{s,Q\varphi^{\ast}})^{\ast}(\lambda)\leq
A_{1}
$$
and
$$
\int_{0}^{\infty}\lambda^{-1}\sum_{Q\in\mathcal{Q}_{\lambda}}|Q|(M_{s,Q\varphi^{\ast}})^{\ast}(\lambda)\leq
A_{1},
$$
there holds
$$
\int_{0}^{\infty}\lambda^{-1}\sum_{Q\in\mathcal{Q}_{\lambda}}|Q|(M_{s,Q\varphi})(\lambda)\leq
A_{2}.
$$
\end{lem}

Note that relation described in  Lemma 8.7 from \cite{D2}  is denoted as
$M_{Q\varphi}\ll(M_{s,Q\varphi^{\ast}})^{\ast}$ (strong domination).

The proof of Lemma 4.8 is based  on some pointwise estimate of
functions $(M_{Q\varphi^{\ast}})^{\ast}$ and
$(M_{s,Q\varphi^{\ast}})^{\ast}.$ This properties we will describe
bellow in Lemma 4.9,4.10.

If $p(\cdot)\in \mathcal{A}^{loc},$  then
$M_{s,Q\varphi}\preceq(M_{s,Q\varphi^{\ast}})^{\ast}(loc)$ for some
$s>1.$ It is not hard to prove that (analogously as  the proof of
Lemma 8.3 from \cite{D2}) uniformly in $Q\in\mathcal{X}^{n}_{loc}$
\begin{equation}
|Q|(M_{s,Q\varphi})\left(\frac{1}{\|\chi_{Q}\|_{p(\cdot)}}\right)\sim1,\,\,\,\,
|Q|(M_{s,Q\varphi}^{\ast})^{\ast}\left(\frac{1}{\|\chi_{Q}\|_{p(\cdot)}}\right)\sim1.
\end{equation}

It is important to investigate for any $Q\in \mathcal{X}^{n}_{loc}$
the function
$$
\alpha_{s}(Q,t)=\frac{(M_{s,Q\varphi})(t)}{(M_{s,Q\varphi^{\ast})^{\ast}}(t)}.
$$

\begin{lem}
Let $p(\cdot)\in \mathcal{A}^{loc}.$ Then uniformly in
$Q\in\mathcal{X}^{n}_{loc}$ and $t>0$
$$
\alpha_{s}(Q,1/\|\chi_{Q}\|_{p(\cdot)})\sim1,\,\,\,\,\alpha_{s}(Q,1)\sim1.
$$
Moreover, there exists $C\geq1$ such that for all
$Q\in\mathcal{X}^{n}_{loc}$
$$
\alpha_{s}(Q,t_{2})\leq
C(\alpha_{s}(Q,t_{1})+1)\,\,\,\,\mbox{for}\,\,0<t_{1}\leq
t_{2}\leq1,
$$

$$
\alpha_{s}(Q,t_{3})\leq
C(\alpha_{s}(Q,t_{4})+1)\,\,\,\,\mbox{for}\,\,1<t_{3}\leq
t_{4}\leq1.
$$
Furthermore, for  all $C_{1},\,C_{2}>0$ there exists $C_{3}\geq1$
such that for all $Q\in\mathcal{X}^{n}_{loc}$
\begin{equation}
t\in\left[C_{1}\min\left\{1,\frac{1}{\|\chi_{Q}\|_{p(\cdot)}}\right\},C_{1}\max\left\{1,\frac{1}{\|\chi_{Q}\|_{p(\cdot)}}\right\}\right]
\,\,\,\Rightarrow\,\,\,\alpha_{s}(Q,t)\leq C_{3}.
\end{equation}
\end{lem}

The proof of analogous statement for nonlocal case (\cite{D2}, Lemma
8.4)is  based on the estimates (4.6)  and some properties (not
depend on $Q$) of convex functions $M_{s,Q\varphi},$
$(M_{s,Q\varphi^{\ast})^{\ast}}.$ This arguments may use in local
variant.

\begin{lem}
Let $p(\cdot)\in \mathcal{A}(loc).$  Then there exists
$b:\,\mathcal{X}^{n}(loc)\rightarrow\mathbb{R}^{\geq}$  and $K>0$
such that
$$
\sup_{\mathcal{Q}\in\mathcal{Y}^{n}(loc)}\sum_{Q\in\mathcal{Q}}|Q|b(Q)+\sup_{Q\in\mathcal{X}^{n}(loc)}|Q|b(Q)<\infty
$$
and for all $Q\in\mathcal{X}^{n}(loc)$ and all $t\geq0$  holds
$$
|Q|(M_{s,Q\varphi^{\ast}})^{\ast}(t)\leq1\,\,\,\Rightarrow\,\,\,(M_{s,Q\varphi})(t)\leq
K(M_{s,Q\varphi^{\ast}})^{\ast}(t)+b(Q).
$$
Moreover, for all $Q\in\mathcal{X}^{n}(loc)$ and all $t\geq1$  there
holds
$$
|Q|(M_{s,Q\varphi^{\ast}})^{\ast}(t)\leq1\,\,\,\Rightarrow\,\,\,(M_{s,Q\varphi})(t)\leq
K(M_{s,Q\varphi^{\ast}})^{\ast}(t).
$$
\end{lem}

The proof may be obtained from general Lemma 4.1 and by using (4.7)
estimate (see \cite{D2}, proof Lemma 8.5).

\begin{lem}
Assume  $M_{s_{2},Q\varphi}\preceq
M_{s_{2},Q\varphi^{\ast}}^{\ast}(loc)$ for some $s_{2}>1$ and $1\leq
s_{1}\leq s_{2}.$ Then uniformly in $Q\in \mathcal{X}^{n}(loc)$ and
$t>0$
$$
\left(\alpha_{s_{2}}(Q,t^{\frac{s_{1}}{s_{2}}})\right)^{\frac{s_{2}}{s_{1}}}\sim
\alpha_{s_{1}}(Q,t).
$$
\end{lem}

The proof of Lemma 4.11 basically based on the Lemma 4.10 and may
proof as analogous  lemma from \cite{D2} (Lemma 8.6).

Let $f$ be locally integrable function. For $q\geq1$ we consider the
local maximal operator  given by
$$
M^{loc}_{q}f(x)=\sup\limits_{Q\ni
x,|Q|\leq1}\left(\frac{1}{|Q|}\int_{Q}|f(y)|^{q}dy\right)^{1/q}.
$$

We define  local dyadic maximal operator $M^{loc}_{q,d}$ witch
restricted supremmum in definition of $M^{loc}_{q}$ by dyadic cubes
(cubes of the form
$Q=2^{-z}((0,1)^{n}+k),\,\,k=(k_{1},...,k_{n})\in\mathbb{Z}^{n},\,\,z\in
\mathbb{N}_{0}).$

For fixed $t\in \mathbb{R}^{n}$  we define also maximal operator
$M^{loc,t}_{q,d}$ whith restricted supremmum in definition of
$M^{loc}_{q,d}$  on the cubes $Q-t,$ where $Q$ dyadic cubes.

Note that there is a constant $C>0$ such that (see \cite{Sa})
\begin{equation}
M^{loc}_{q}f(x)\leq C\int_{[-4,4]^{n}}M^{loc,t}_{q,d}f(x)dt.
\end{equation}

The main step to proof Theorem 2.2 (as in proof of original Theorem
2.1) is following Theorem.
\begin{thm}
Let $p(\cdot)\in \mathcal{A}_{loc}.$  Then there exists $q>1$ such
that $M_{q}^{loc}$  is continuous on $L^{p(\cdot)}(\mathbb{R}^{n}).$
\end{thm}

Note that using (4.8) estimate it is sufficiently to proof Theorem
4.12 for operator $M^{loc}_{q,d}.$

It is suffices to show that there exists $A>0$ such that for all
$f\in L^{p(\cdot)}(\mathbb{R}^{n})$
$$
\int_{\mathbb{R}^{n}}|f(x)|^{p(x)}dx\leq1\,\,\,\Rightarrow\,\,\int_{\mathbb{R}^{n}}|M^{loc}_{q,d}f(x)|^{p(x)}dx\leq
A.
$$

  For $\lambda>0$ define functions
  $$
  f_{0,\lambda}=f\chi_{\{|f|\leq\lambda}\},\,\,\,f_{0,\lambda}=f\chi_{\{|f|>\lambda\}}.
  $$
  Then
  $$
  \{M^{loc}_{q,d}f>\lambda\}\subset\{M^{loc}_{q,d}f_{0,\lambda}>\lambda/2\}\cup\{M^{loc}_{q,d}f_{1,\lambda}>\lambda/2\}.
  $$

This implies
$$
\int_{\mathbb{R}^{n}}|M^{loc}_{q,d}f(x)|^{p(x)}dx=\int_{0}^{\infty}\int_{\mathbb{R}^{n}}p(x)\lambda^{p(x)-1}\chi_{\{M^{loc}_{q,d}f>\lambda\}}dxd\lambda
$$
$$
\leq
C\sum_{j=1}^{2}\int_{0}^{\infty}\lambda^{-1}\int_{\mathbb{R}^{n}}\lambda^{p(x)}\chi_{\{M^{loc}_{q,d}f_{j,\lambda}>\lambda/2\}}dxd\lambda.
$$
For $\lambda>0$ let $\mathcal{Q}_{0,\lambda}$ be the decomposition
of $\{M^{loc}_{q,d}f_{0,\lambda}>\lambda/2\}$ into maximal dyadic
cubes. Then for all $Q\in\mathcal{Q}_{0,\lambda}$ there holds
(uniformly in $Q$)
$$
M_{q,Q}f_{0,\lambda}\sim \lambda
$$
and we have
$$
\int_{0}^{\infty}\lambda^{-1}\int_{\mathbb{R}^{n}}\lambda^{p(x)}\chi_{\{M^{loc}_{q,d}f_{0,\lambda}>\lambda/2\}}dxd\lambda\leq
C\int_{0}^{\infty}\lambda^{-1}\sum_{Q\in\mathcal{Q}_{0,\lambda}}|Q|(M_{Q\varphi})(\lambda)d\lambda.
$$

Denote
$f_{1,\lambda}^{k}=\chi_{(0,1)^{n}+k}f_{1,\lambda},\,k\in\mathbb{Z}^{n}.$
Note that if $x\in(0,1)^{n}+k$ then
$$
M^{loc}_{q,d}f_{1,\lambda}(x)=M^{loc}_{q,d}f_{1,\lambda}^{k}(x)
$$
and
$$
\{M^{loc}_{q,d}f_{1,\lambda}>\lambda/2\}=\cup_{k\in\mathbb{Z}^{n}}\{M^{loc}_{q,d}f_{1,\lambda}^{k}>\lambda/2\}.
$$

We have
$$
\int_{0}^{\infty}\lambda^{-1}\int_{\mathbb{R}^{n}}\lambda^{p(x)}\chi_{\{M^{loc}_{q,d}f_{1,\lambda}>\lambda/2\}}dxd\lambda
$$
$$
=\int_{0}^{\infty}\lambda^{-1}\sum_{k\in\mathbb{Z}^{n}}\int_{(0,1)^{n}+k}\lambda^{p(x)}\chi_{\{M^{loc}_{q,d}f_{1,\lambda}^{k}>\lambda/2\}}dxd\lambda.
$$

Define $m_{k}=2\int_{(0,1)^{n}+k}|f(x)|dx,$  we have

$$
\int_{0}^{\infty}\lambda^{-1}\int_{(0,1)^{n}+k}\lambda^{p(x)}\chi_{\{M^{loc}_{q,d}f_{1,\lambda}^{k}>\lambda/2\}}dxd\lambda
$$

$$
=\int_{0}^{m_{k}}\lambda^{-1}\int_{(0,1)^{n}+k}\lambda^{p(x)}\chi_{\{M^{loc}_{q,d}f_{1,\lambda}^{k}>\lambda/2\}}dxd\lambda
$$
$$
+
\int_{m_{k}}^{\infty}\lambda^{-1}\int_{(0,1)^{n}+k}\lambda^{p(x)}\chi_{\{M^{loc}_{q,d}f_{1,\lambda}^{k}>\lambda/2\}}dxd\lambda.
$$

Note that
$$
\int_{0}^{m_{k}}\lambda^{-1}\int_{(0,1)^{n}+k}\lambda^{p(x)}\chi_{\{M^{loc}_{q,d}f_{1,\lambda}^{k}>\lambda/2\}}dxd\lambda\leq
C\int_{(0,1)^{n}+k}\left(\int_{(0,1)^{n}+k}|f(t)|dt\right)^{p(x)}dx.
$$

Let $\mathcal{Q}_{1,\lambda}^{k}$ be the decomposition of
$\{M^{loc}_{q,d}f_{1,\lambda}^{k}>\lambda/2\}$ into maximal dyadic
cubes. Then for all $Q\in\mathcal{Q}_{1,\lambda}^{k}$ there holds
$$
M_{q,Q}f_{1,\lambda}^{k}=M_{q,Q}f_{1,\lambda}\sim\lambda.
$$
Define
$\mathcal{Q}_{1,\lambda}=\cup_{k\in\mathbb{Z}^{n}}\mathcal{Q}_{1,\lambda}^{k}.$
Then we have
$$
\int_{0}^{\infty}\lambda^{-1}\int_{\mathbb{R}^{n}}\lambda^{p(x)}\chi_{\{M^{loc}_{q,d}f_{1,\lambda}>\lambda/2\}}dxd\lambda
$$
$$
\leq
C\sum_{k\in\mathbb{Z}^{n}}\int_{(0,1)^{n}+k}\left(\int_{(0,1)^{n}+k}|f(t)|dt\right)^{p(x)}dx+
\int_{0}^{\infty}\lambda^{-1}\sum_{Q\in\mathcal{Q}_{1,\lambda}}|Q|(M_{Q\varphi})(\lambda)d\lambda.
$$

For first term we have
$$
\sum_{k\in\mathbb{Z}^{n}}\int_{(0,1)^{n}+k}\left(\int_{(0,1)^{n}+k}|f(t)|dt\right)^{p(x)}dx\leq
C.
$$
The second term
$\int_{0}^{\infty}\lambda^{-1}\sum_{Q\in\mathcal{Q}_{1,\lambda}}|Q|(M_{Q\varphi})(\lambda)d\lambda$
can be estimated in the same way as in the Theorem 6.2 from
\cite{D2}.


\end{document}